\documentclass[reqno,12pt]{amsart}

\usepackage{psfig, amsmath, amsfonts, amssymb, amsthm, amscd, graphics,epsf,psfrag}

\numberwithin{equation}{section}

\makeatletter
\def\captionfont@{\footnotesize}
\def\captionheadfont@{\scshape}

\long\def\@makecaption#1#2{%
  \vspace{2mm}
  \setbox\@tempboxa\vbox{\color@setgroup
    \advance\hsize-6pc\noindent
    \captionfont@\captionheadfont@#1\@xp\@ifnotempty\@xp
        {\@cdr#2\@nil}{.\captionfont@\upshape\enspace#2}%
    \unskip\kern-6pc\par
    \global\setbox\@ne\lastbox\color@endgroup}%
  \ifhbox\@ne % the normal case
    \setbox\@ne\hbox{\unhbox\@ne\unskip\unskip\unpenalty\unkern}%
  \fi
  \ifdim\wd\@tempboxa=\z@ % this means caption will fit on one line
    \setbox\@ne\hbox to\columnwidth{\hss\kern-6pc\box\@ne\hss}%
  \else % tempboxa contained more than one line
    \setbox\@ne\vbox{\unvbox\@tempboxa\parskip\z@skip
        \noindent\unhbox\@ne\advance\hsize-6pc\par}%
\fi
  \ifnum\@tempcnta<64 % if the float IS a figure...
    \addvspace\abovecaptionskip
    \moveright 3pc\box\@ne
  \else % if the float IS NOT a figure...
    \moveright 3pc\box\@ne
    \nobreak
    \vskip\belowcaptionskip 
  \fi
\relax
}
\makeatother
\def\writefig#1 #2 #3 {\rlap{\kern #1 truecm
\raise #2 truecm \hbox{#3}}}

\psfigurepath{.:./pictures}

\DeclareMathSymbol{\leqslant}{\mathalpha}{AMSa}{"36} % nicer `smaller or equal'
\DeclareMathSymbol{\geqslant}{\mathalpha}{AMSa}{"3E} % nicer `larger or equal'
\DeclareMathSymbol{\eset}{\mathalpha}{AMSb}{"3F}     % nicer `emptyset'
\renewcommand{\leq}{\;\leqslant\;}                   % redef. of < or =
\renewcommand{\geq}{\;\geqslant\;}                   % redef. of > or =
\newtheorem{pro}{Proposition}[section]
\newtheorem{thm}{Theorem}[section]

\newtheorem{defi}{Definition}[section]

\newcommand{\cA}{\ensuremath{\mathcal A}}
\newcommand{\cB}{\ensuremath{\mathcal B}}
\newcommand{\cC}{\ensuremath{\mathcal C}}
\newcommand{\cD}{\ensuremath{\mathcal D}}
\newcommand{\cE}{\ensuremath{\mathcal E}}

\newcommand{\cT}{\ensuremath{\mathcal T}}

\newcommand{\cZ}{\ensuremath{\mathcal Z}}

\newcommand{\frC}{\ensuremath{\mathfrak C}}

\newcommand{\bbB}{{\ensuremath{\mathbb B}} }

\newcommand{\bbE}{{\ensuremath{\mathbb E}} }

\newcommand{\bbP}{{\ensuremath{\mathbb P}} }
\newcommand{\bbQ}{{\ensuremath{\mathbb Q}} }

\newcommand{\bbZ}{{\ensuremath{\mathbb Z}} }

\newcommand{\ga}{\alpha}
\newcommand{\gb}{\beta}
\newcommand{\gga}{\gamma}            % \gg already exists...

\newcommand{\gG}{\Gamma}

\newcommand{\go}{\omega}
\newcommand{\gO}{\Omega}
\newcommand{\gl}{\lambda}
\newcommand{\gL}{\Lambda}
\newcommand{\gs}{\sigma}

\newcommand{\Perc}{\Phi}           % Percolation measure
\def\1{\ifmmode {1\hskip -3pt \rm{I}} \else {\hbox {$1\hskip -3pt \rm{I}$}}\fi}

\newcommand{\lra}{\leftrightarrow}

\newcommand{\FKm}[2]{\Perc^{#1}_{#2}}

\newcommand{\Ew}[1]{\bbE^{{\rm w}}_{#1}}
\newcommand{\Ef}[1]{\bbE^{{\rm f}}_{#1}}

\title[Translation invariant Gibbs states]{Translation invariant Gibbs states\\ 
for the Ising model}

\author{T. Bodineau}
\address{D{\'e}partement de math{\'e}matiques, Universit{\'e} Paris 7,
case 7012, 2 place Jussieu, Paris 75251, France}
\email{bodineau@gauss.math.jussieu.fr}
\thanks{}

\date{\today}

\begin{document}

\maketitle

\begin{abstract}
We prove that all the translation invariant Gibbs states of the Ising
model are a linear combination of the pure phases
$\mu^+_\gb,\mu^-_\gb$ for any $\gb \not = \gb_c$.
This implies that the average magnetization is continuous for $\gb
>\gb_c$. Furthermore, combined with previous results on the slab percolation
threshold \cite{Bo2} this shows the validity of Pisztora's coarse
graining  \cite{pisztora} up to the critical temperature.
\end{abstract}

\section{Introduction}

The set of Gibbs measures associated to the Ising model is a simplex
(see \cite{Georgii}) and the complete characterization of the extremal
measures at any inverse temperature $\gb = 1/T$ remains an important 
issue. 
The most basic states are the two pure phases $\mu^+_\gb,\mu^-_\gb$
which are obtained as the thermodynamic limit of the finite Gibbs
measures with boundary conditions uniformly equal to $1$ or $-1$.
In the phase transition regime ($\gb > \gb_c$), these two Gibbs 
states are distinct and translation invariant.
An important result by Aizenman and Higuchi \cite{Aizenman,Hig} (see 
also \cite{GeoHig}) asserts that for the two dimensional nearest neighbor
Ising model these are the only two extremal Gibbs measures and that
any other Gibbs measure on $\{\pm 1\}^{\bbZ^2}$ belongs to 
$[\mu^+_\gb,\mu^-_\gb]$, i.e. is a linear combination of 
$\mu^+_\gb,\mu^-_\gb$.
In higher dimensions Dobrushin \cite{D} proved the existence of other 
extremal invariant measures. They arise from well chosen mixed
boundary conditions which create a rigid interface separating the
system into two regions. Thus, contrary to the previous pure phases,
the Dobrushin states are non-translation invariant. We refer the
reader to the survey by Dobrushin, Shlosman \cite{DS} for a detailed
account on these states.

\medskip

In this paper we are going to focus on the translation invariant 
Gibbs states in the phase transition regime
and prove that they belong to $[\mu^+_\gb,\mu^-_\gb]$.
This problem has a long history and has essentially already been 
solved, with the exception of one detail which we will now tie up.

Two strategies have been devised to tackle the problem.
The first one, implemented by Gallavotti and Miracle-Sol\'e \cite{GS},
is a constructive method based on Peierls estimates.
They proved that for any $\gb$ large enough the 
set of translation invariant Gibbs states is $[\mu^+_\gb,\mu^-_\gb]$.
This result was generalized in \cite{BMP} to the Ising model with 
Kac interactions for any $\gb >1$ as soon as the interaction range
is large enough. 
A completely different approach relying on ferromagnetic inequalities
was introduced by Lebowitz \cite{Leb1} and generalized to the framework
of FK percolation by Grimmett \cite{G2}.
The key argument is to relate the differentiability of the pressure
wrt $\gb$ and the characterization of the translation invariant 
Gibbs states.
As the pressure is a convex function, it is differentiable for all
$\gb$, except possibly for an at most countable set of inverse temperatures
$\cB \subset [\gb_c, \infty[$.
For the Ising model, $\cB$ is conjectured to be empty, although
the previous method does not provide any explicit control on $\cB$.
We stress the fact that the non differentiability of the pressure
has other implications, namely that for any inverse
temperature in $\cB$, the average magnetization would be discontinuous;
and that the number of pure phases would be uncountable (see \cite{BL}).

We will show that for any $\gb>\gb_c$ there is a unique infinite
volume FK measure. Several consequences can be drawn from this
by using previous results in \cite{G2,Leb1}: the set of translation invariant 
Gibbs states is  $[\mu^+_\gb,\mu^-_\gb]$, the average magnetization
is continuous in $]\gb_c, \infty[$.
Finally, combining this statement with the characterization of the
slab percolation threshold in \cite{Bo2}, we deduce that  Pisztora's
coarse graining is valid up to the critical temperature.
All these facts are summarized in Subsection \ref{subsec: Results}.
Our method is restricted to $\gb >\gb_c$. However,
it is widely believed that the phase transition of the Ising model 
is of second order and thus similar results should also hold at
$\gb_c$.

\section{Notation and Results}

	\subsection{The Ising model}
	\label{subsec: Ising}

We consider the Ising model on $\bbZ^d$ with finite range interactions
and  spins $\{ \gs_i \}_{ i \in \bbZ^d}$
taking values  $\pm 1$. 
Let $\gs_{\gL} \in \{ \pm 1 \}^{\gL}$ be the spin configuration 
restricted to $\gL \subset \bbZ^d$.
The Hamiltonian associated to $\gs_{\gL}$ with boundary conditions 
$\gs_{\gL^c}$ is defined by
\begin{eqnarray*}
H( \gs_{\gL}  \, | \, \gs_{\gL^c} ) = 
- {1 \over 2} \sum_{ i,j \in \gL} J(i-j) \gs_i \gs_j
- \sum_{ i \in \gL, j \in \gL^c} J(i-j) \gs_i \gs_j,
\end{eqnarray*}
where the couplings $J(i-j)$ are ferromagnetic and equal to 0 
for $\| i -j \| \geq R$ ($R$ will be referred to as the range
of the interaction).

The Gibbs measure in $\gL$ at inverse temperature $\gb > 0$ is
defined by
\begin{eqnarray*}
\label{Gibbs measure}
\mu_{\gb, \gL}^{\gs_{\gL^c}} ( \gs_{\gL} ) =
{1 \over Z_{\gb,\gL}^{\gs_{\gL^c}} }
\exp \big( - \gb  H( \gs_{\gL}  \, | \, \gs_{\gL^c} ) \big), 
\end{eqnarray*}
where the partition function $Z_{\gb,\gL}^{\gs_{\gL^c}}$ is the normalizing 
factor. 
The boundary conditions act as  boundary fields, therefore more general
values of the boundary conditions can be used. For any $h>0$, let us denote 
by $\mu_{\gb, \gL}^{h} $ the Gibbs measure with boundary magnetic field
$h$, i.e. with Hamiltonian
\begin{eqnarray*}
H_h ( \gs_{\gL} ) = 
- {1 \over 2} \sum_{i,j \in \gL} J(i-j) \gs_i \gs_j
-  h \sum_{i \in \gL, j \in \gL^c} J(i-j) \gs_i \, .
\end{eqnarray*}

\medskip

The phase transition is characterized by symmetry breaking for any
$\gb$ larger than  the inverse critical temperature $\gb_c$ defined by
\begin{equation*}
\gb_c = \inf \{ \gb>0, \qquad 
\lim_{N\to\infty} \mu_{\gb,\gL_N}^+ (\gs_0) >0 \} \, . 
\end{equation*}

	\subsection{The random cluster measure}
	\label{subsec: FK}

The random cluster measure was originally introduced 
by Fortuin and Kasteleyn \cite{FK} (see also \cite{ES,G2})
and it can be understood as an alternative representation 
of the Ising model (or more generally of the $q$-Potts model).
This representation will be referred to as the FK representation.

Let $\bbE$ be the set of bonds, i.e. of pairs $(i,j)$ in $\bbZ^d$ such
that $J(i-j)>0$.
For any subset $\gL$ of $\bbZ^d$ we consider two sets of bonds
\begin{equation}
\label{eq: edges}
\begin{cases}
\Ew{\gL} = \{ (i,j) \in \bbE, \quad i \in \gL, j \in \bbZ^d \} \, ,\\
\Ef{\gL} = \{ (i,j) \in \bbE, \quad i,j \in \gL \} \, .
\end{cases}
\end{equation}

The set $\Omega = \{ 0,1\}^{\bbE}$  is the state space 
for the dependent percolation measures.
Given $\go\in\Omega$ and a bond $b=(i,j) \in\bbE$, we say  that $b$ is open 
if $\go_b=1$. Two sites of $\bbZ^d$ are said
to be connected if one can be reached from another via a chain of open bonds.
Thus, each $\go\in\Omega$ splits  $\bbZ^d$ into the disjoint union of maximal
connected components, which are called the open clusters of $\Omega$. Given a
finite subset $B\subset\bbZ^d$  we use $c_B (\go )$ to denote the number of
different open finite 
clusters of $\go$ which have a non-empty intersection
with $B$.

For any $\gL \subset \bbZ^d$ we define the random cluster measure 
on the bond configurations $\omega\in \Omega_\gL  =  \{0
,1\}^{\Ef{\gL}}$.
The boundary conditions are specified by a frozen percolation configuration  
$\pi \in \Omega_\gL^c = \Omega \setminus \Omega_\gL $. 
Using the shortcut  $c^\pi_\gL (\go ) =c_{\gL} (\go\vee \pi )$ for 
the joint configuration $\go \vee \pi \in \bbE$, we define the finite 
volume random cluster measure $\FKm{\pi}{\gb,\gL}$ on $\Omega_\gL$  with the 
boundary conditions $\pi$ as:
\begin{equation}
\label{FKm}
\FKm{\pi}{\gb,\gL}\left(\go \right)~ = ~\frac1{Z^{\gb ,\pi}_\gL}
\left( \prod_{b \in \Ef{\gL} } \big( 1-p_b \big)^{1-\go_b} \; p_b^{\go_b} 
\right) \, 2^{c^\pi_\gL (\go)}\,,
\end{equation}
where the bond intensities are such that $p_{(i,j)} = 1-\exp(-2\gb J(i-j))$.
We will sometimes use the same notation for the FK measure on $\Ew{\gL}$, 
in which case we will state it explicitly.

\medskip

The measures $\FKm{\pi}{\gb,\gL}$ are FKG partially ordered with respect to
the lexicographical order of the boundary condition $\pi$. Thus, the extremal
ones correspond to the free ($\pi\equiv 0$) and wired ($\pi\equiv 1$) boundary
conditions and are denoted as  $\FKm{{\rm f}}{\gb,\gL}$ and $\FKm{{\rm
w}}{\gb,\gL}$  respectively. 
The  corresponding infinite volume limits $\FKm{{\rm f}}{\gb}$ and 
$\FKm{{\rm w}}{\gb}$ always exist.

The phase transition of the random cluster model is characterized by 
the occurrence of percolation 
\begin{equation}
\label{eq: transition}
\forall \gb > \gb_c, \qquad 
\lim_{N\to\infty}\FKm{{\rm w}}{\gb, \gL_N}\left( 0 \lra
\gL_N^c\right)~=~ \FKm{ {\rm w}}{\gb} \left( 0\lra\infty \right)
> 0.
\end{equation}
%The correspondence with the Ising model implies that $p_c = 1-\exp(-2\gb_c)$.

		\subsection{Results and consequences}
		\label{subsec: Results}

Our main result is  
\begin{thm}
\label{thm: main}
In the case of Ising model  for any $\gb \not = \gb_c$
\begin{eqnarray}
\label{eq: egalite}
\FKm{{\rm f}}{\gb} \big( \{ 0 \lra \infty \} \big)
= \FKm{{\rm w}}{\gb} \big( \{ 0 \lra \infty \} \big) \, .
\end{eqnarray}
\end{thm}

The proof is postponed to Subsection \ref{subsec: proof main thm} and we first draw 
some consequences from this Theorem.

\medskip

\noindent
$\bullet$ {\bf Continuity of the average magnetization.}

Grimmett proved in  \cite{G2} (Theorem 5.2) that the function
$\gb \to \FKm{{\rm w}}{\gb} \big( \{ 0 \lra \infty \} \big)$ is 
right continuous in $[0,1]$ 
and $\gb \to \FKm{{\rm f}}{\gb} \big( \{ 0 \lra \infty \} \big)$ is
left continuous in $[0,\infty[ \setminus \{\gb_c\}$. 
Therefore Theorem \ref{thm: main} implies that the average magnetization 
\begin{eqnarray}
\label{eq: averaged magnetization}
\mu^+_\gb (\gs_0) 
= \FKm{{\rm w}}{\gb} \big( \{ 0 \lra \infty \} \big)
%, \qquad \text{with} \quad p = 1 - \exp(-2\gb)  
\end{eqnarray}
is a continuous function of $\gb$ except possibly at $\gb_c$.

\medskip

\noindent
$\bullet$ {\bf Translation invariant states.}

According to Theorem 5.3 (b) in  \cite{G2}, equality \eqref{eq: egalite} 
implies that there exists only one random cluster measure. 
This means that $\FKm{{\rm w}}{\gb} = \FKm{{\rm f}}{\gb}$ for $\gb \not = \gb_c$.

Alternatively for the spin counterpart, Lebowitz proved in \cite{Leb1}
(Theorem 3 and remark (iii) page 472) 
that the continuity of the average magnetization implies the existence of only two
extremal invariant states, i.e. that for $\gb > \gb_c$ all the
translation invariant Gibbs 
states are of the form $\gl \mu^+_\gb + (1 -\gl) \mu^-_\gb$ for some $\gl \in [0,1]$.

\medskip

\noindent
$\bullet$ {\bf Pisztora's coarse graining.}

A description of the Ising model close to the critical temperature
requires a renormalization procedure in order to deal with the
diverging correlation length.
A crucial tool for implementing this is the Pisztora's coarse graining
\cite{pisztora} which provides an accurate description of the typical
configurations of the Ising model (and more generally of the
$q$-Potts model) in terms of the FK representation.
This renormalization scheme is at the core of many works on the Ising
model and in particular it was essential for the analysis of phase 
coexistence (see \cite{Ce,CePi,Bo,BIV}).

The main features of the coarse graining will be recalled in Subsection 
\ref{subsec: Renormalization}.
Nevertheless, we stress that its implementation is based upon two
hypothesis:
\begin{enumerate}
\item The inverse temperature $\gb$ should be above the slab
percolation threshold (see \cite{pisztora}).
\item The uniqueness of the FK measure, i.e. 
$\FKm{{\rm f}}{\gb}= \FKm{{\rm w}}{\gb}$.
\end{enumerate}

The first assumption was proved to hold for the Ising model as soon as
$\gb > \gb_c$ \cite{Bo2} and as a consequence of Theorem 
\ref{thm: main}, the second is also valid for $\gb > \gb_c$.
Thus for the Ising model, Pisztora's coarse graining applies in the
whole of the phase transition regime and from \cite{CePi} the Wulff 
construction in dimension $d \geq 3$ is valid up to the critical temperature.

	\section{Proof of Theorem \ref{thm: main}}

Let us briefly comment on the structure of the proof.
It is well known that the wired measure $\FKm{{\rm w}}{\gb}$ 
dominates the free measure $\FKm{{\rm f}}{\gb}$ in the FKG sense
thus the core of the proof is to prove the reverse inequality.
The first step is to show that $\FKm{{\rm f}}{\gb}$ dominates 
the FK counterpart of the finite volume Gibbs measure $\mu^h_{\gb,
\gL}$ for some value of $h>0$ and independently of $\gL$.
This is achieved by introducing intermediate random variables $Z$
(Subsection \ref{subsec: Free boundary conditions}) and 
$\widehat Z$ (Subsection \ref{subsec: wired}) which can be compared
thanks to a coupling (Subsection \ref{subsec: coupling}). 
We then rely on a result by Lebowitz \cite{Leb2} 
and Messager, Miracle Sole, Pfister \cite{MMP} which ensures
that $\mu^h_{\gb, \gL}$ converges to $\mu^+_{\gb}$ in the
thermodynamic limit.
From this, we deduce that $\FKm{{\rm f}}{\gb}$ dominates 
$\FKm{{\rm w}}{\gb}$ in the FKG sense (Subsection \ref{subsec: proof main thm}).

		\subsection{Renormalization}
		\label{subsec: Renormalization}

We recall the salient features of Pisztora's coarse graining 
and refer to the original paper \cite{pisztora} for the details.
The reference scale for the coarse graining is  an integer $K$ which
will be chosen large enough.
The space $\bbZ^d$ is partitioned into blocks of side length $K$
\begin{eqnarray*}
\forall x \in K\, \bbZ^d, \qquad 
\bbB_K (x) = x + \left\{ - \frac{K}{2}+1, \dots, \frac{K}{2} \right\}^d
\, .
\end{eqnarray*}

First of all we shall set up the notion of {\it good} block on the $K$-scale which characterizes a local equilibrium in a pure phase.

\begin{defi}
A block $\bbB_K (x)$ is said to be {\it good} with respect to the bond configuration 
$\go \in \gO$ if the following events are satisfied
\begin{enumerate}
\item There exists a crossing cluster ${\bf C}^*$ in $\bbB_K (x)$ connected to all 
the faces of the inner vertex boundary of $\bbB_K (x)$. 
\item Any FK-connected cluster in  $\bbB_K (x)$ of diameter larger than $\sqrt{K}/10$
is contained in ${\bf C}^*$.
\item There are crossing clusters in each block $\big( \bbB_{\sqrt{K}} (x \pm 
\frac{K}{2} \vec{e}_i) \big)_{1 \leq i \leq d}$, where $\big( \vec{e}_i \big)_{1 \leq i \leq d}$ 
are the unit vectors (see (4.2) in \cite{pisztora}).
\item There is at least a closed bond in $\bbB_{K^{1/2d}} (x)$.
\end{enumerate}
\end{defi}

The important fact which can be deduced from (1,2,3) is that the crossing clusters
in two neighboring good blocks are connected. Thus a connected cluster of good blocks at scale $K$ induces also the occurrence of  a connected cluster at the microscopic level.

To each block $\bbB_K (x)$, we associate a  coarse grained variable $u_K(x)$ equal to 1
if this is a good block or 0 otherwise.
Fundamental techniques developed by Pizstora (see (4.15) in \cite{pisztora}) imply that
a block is good with high probability conditionally to the states of
its neighboring blocks.
For any $\gb > \gb_c$, there is $K_0$ large enough such that for all scales 
$K\geq K_0$ one can find a constant $C > 0$ (depending on $K,\gb$) such that 
\begin{eqnarray} 
\label{eq: Peierls} 
\FKm{{\rm f}}{\gb}  
\left( u_K(x) = 0 \; \Big | \; u_K(y) = \eta_y, \quad y \not = x  \right) 
\leq  \exp ( - C) \, ,
\end{eqnarray} 
this bound holds uniformly over the values $\eta_j \in \{0,1\}$ of the neighboring
blocks. 
Furthermore, the constant $C$ diverges as $K$ tends to infinity.
The previous estimate was originally derived beyond the slab percolation threshold.
The latter has been proved to coincide with the critical temperature in the case
of the Ising model \cite{Bo2}.

A last feature of Pisztora's coarse graining is a control of the density of the 
crossing cluster in each good block. Under the assumption that \eqref{eq: egalite}
holds, one can prove that with high probability, the density of the crossing cluster 
in each block is close to the one of the infinite cluster.
Thus,  one of the goals of this paper is to prove that the complete renormalization scheme is valid up to the critical temperature.
Throughout the paper, we will use only the estimate \eqref{eq: Peierls} and not
the full Pisztora's coarse graining which includes as well the control on the
density.

\medskip

\begin{figure}[h]
\begin{center}
\leavevmode
\epsfysize = 5 cm
\psfrag{L}[B]{$L$}
\psfrag{N}[Br]{$N$}
\psfrag{D}[Bl]{$\gL_N$}
\psfrag{X}[l]{$\Xi_{N,L}$}
\psfrag{T}{$T_L(x)$}
\epsfbox{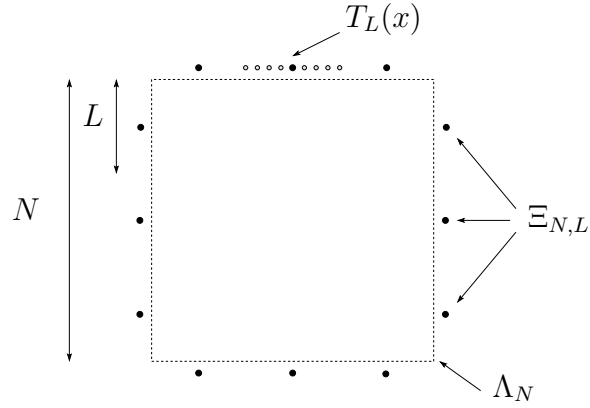}
\end{center}
\caption{The figure corresponds to the nearest neighbor Ising model.
The scales are not accurate and one should imagine $1 \ll K \ll L \ll N$.
The set $\gL_N$ is depicted in dashed lines. The subset $\Xi_{N,L}$ is the
union of the black dots which all  belong to $\partial \gL_N$. Only one set $T_L(x)$
has been depicted at the top.}
\label{fig: figure1}
\end{figure}

For $N = n \frac{K}{2}$, we define 
\begin{eqnarray}
\label{eq: set notation}
\gL_N = \{ -N +1 ,\dots, N\}^d, \qquad
\partial \gL_N = \{ j \in \gL_N^c \; |  \quad \exists i \in \gL_N, \ 
J(i-j)>0 \} \, .
\end{eqnarray}
The set $\partial \gL_N$ is the boundary of $\gL_N$. 
It will be partitioned into $(d-1)$-dimensional 
slabs of side length $L = \ell K$ (for some appropriate choice of
$n$ and $\ell$). 
More precisely if $R$ denotes the range of the interaction, we define  the slab
$$T_L = \{0,\dots,R\} \times \{-L/2 + 1,\dots, L/2\}^{d-1}$$ 
and  $\Xi_{N,L}$ a subset of $\partial \gL_N$ 
such that $\partial \gL_N$
can be covered by non intersecting slabs with centers
in $\Xi_{N,L}$
\begin{eqnarray}
\label{eq: partition}
\partial \gL_N = \bigcup_{x \in \Xi_{N,L}} T_L(x) \, ,
\end{eqnarray}
where $T_L(x)$ denotes the slab centered at site $x$ and
deduced from $T_L$ by  rotation and translation
(see figure \ref{fig: figure1}).

		\subsection{Free boundary conditions}
		\label{subsec: Free boundary conditions}

\begin{figure}[h]
\begin{center}
\leavevmode
\epsfysize = 5 cm
\psfrag{L}[Br]{$\gL_{N + 3K/2}^c$}
\psfrag{x}[Br]{$x$}
\psfrag{y}[Bl]{$y$}
\psfrag{G}[Tr]{Path of good Blocks}
\psfrag{T}[r]{$T_L(x)$}
\psfrag{I}[l]{$\infty$}
\epsfbox{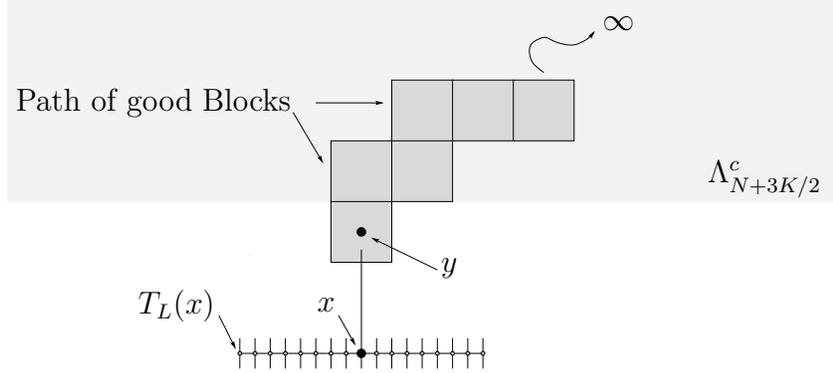}
\end{center}
\caption{The event $Z_x = 1$ is depicted (the scales are not accurate).
The black lines are
the open bonds attached to $T_L(x)$. The block $\bbB_K(y)$ is good and
connected to infinity by a path of good blocks included in $\gL_{N +
3K/2}^c$ (represented by the light gray region).}
\label{fig: Z}
\end{figure}

We define new random variables indexed by the set $\Xi_{N,L}$
introduced in \eqref{eq: partition}.
\begin{defi}
\label{def: Z}
The collection $(Z_x)_{x \in \Xi_{N,L}}$ depends on the bond
configurations in $\bbE \setminus \Ef{\gL_N}$.
For any $x$ in $\Xi_{N,L}$, we declare that $Z_x=1$ if 
the three following events are satisfied (see figure \ref{fig: Z})
\begin{enumerate}
\item All the bonds in $\bbE \setminus \Ef{\gL_N}$ intersecting $T_L(x)$ are open.
\item If $\vec{n}$ denotes the outward normal to $\gL_{N+1}$ at $x$
then the $3K/4$ edges $\big\{ \big(x + i \vec{n}, x + (i+1) \vec{n}
\big) \big\}_{0 \leq i \leq 3K/4}$ are open. 
Let $y$ be the site $x+  K \vec{n}$. Then $\bbB_K(y)$ is a
good block, i.e. $u_K(y) = 1$.
\item 
The block $\bbB_K(y)$ is connected to infinity by an
open path of good blocks included in
$\gL_{N+ 3K/2}^c$. 
\end{enumerate}
If one of the events is not satisfied, then $Z_x=0$.

\medskip

Let $\bbQ$ be the image measure on $\{0,1\}^{\Xi_{N,L}}$ of 
$\FKm{{\rm f}}{\gb}$ by the application $\go \to \{ Z_x(\go) \}_{x \in \Xi_{N,L}}$.
\end{defi}

It is convenient to order the sites of $\Xi_{N,L}$ wrt the lexicographic 
order and to index the random variables by $\{Z_k\}_{k \leq M}$,
where $M$ is the cardinality of  $\Xi_{N,L}$. The $k^{th}$ element
$x_k$ of $\Xi_{N,L}$ is associated to $Z_k = Z_{x_k}$.

\medskip

We will associate to a given sequence $\{Z_k\}_{k \leq M}$ a random cluster measure in $\Ef{\gL_N}$ with boundary conditions which will be wired in the regions where $Z_k =1$
and free otherwise.
More precisely, $\partial \gL_N$ is split into two regions
\begin{equation*}
\partial^{{\rm f}} \gL_N =\bigcup_{k \ \text{such that} \ Z_k = 0} T_L (x_k),
\qquad    
\partial^{{\rm w}} \gL_N = \bigcup_{k \ \text{such that} \ Z_k = 1} T_L (x_k) \, .
\end{equation*}
We set 
\begin{equation}
\label{eq: pi}
\forall (i,j) \in \Ew{\gL_N} \setminus \Ef{\gL_N}, 
\qquad
\pi^Z_{(i,j)} =
\begin{cases}
0, \quad \text{if} \ \  i \in \partial^{{\rm f}} \gL_N, \  j \in \gL_N \,  ,  \\
1, \quad \text{if} \ \  i \in \partial^{{\rm w}} \gL_N, \  j \in \gL_N \,  .
\end{cases}
\end{equation}
Outside $\Ew{\gL_N}$ the boundary conditions will be wired and we
set $\pi^Z_b = 1$ for $b$ in $\bbE \setminus \Ew{\gL_{N}}$.
Finally, let us introduce for the FK measure in $\Ef{\gL_N}$ with
boundary conditions $\pi^Z$ 
\begin{eqnarray}
\label{eq: FK Z}
\forall Z \in \{0,1\}^{\Xi_{N,L}}, \qquad 
\Psi(Z) = 
\FKm{\pi^Z}{\gb,\gL_N} \big( 0 \lra \partial^{{\rm w}} \gL_N \big) \, .
\end{eqnarray}
If $\partial^{{\rm w}} \gL_N$ is empty then $\Psi (Z) = 0$.

By construction, to any bond configuration $\go$ outside 
$\Ef{\gL_N}$, one can associate a collection $\{Z_k (\go)\}$ and 
a bond configuration $\pi^{Z(\go)}$. 
Almost surely wrt $\FKm{{\rm f}}{\gb}$, the infinite cluster is unique
for any $\gb > \gb_c$ \cite{BK} and all the sites $x_k$ such that $Z_k =1$
belong to the same cluster. 
Thus the following FKG domination holds
\begin{eqnarray*}
\FKm{\go}{\gb,\gL_N} \succ \FKm{\pi^{Z(\go)}}{\gb,\gL_N}, \qquad   
\FKm{{\rm f}}{\gb} \   a.s.
\end{eqnarray*}
As the event $\{ 0 \lra \infty\}$ is increasing, we get
\begin{eqnarray}
\label{eq: omega > Z}
\FKm{{\rm f}}{\gb} \big( 0 \lra \infty \big)
\geq 
\bbQ \left(\Psi \big( Z \big) \right) \, . 
\end{eqnarray}

\medskip

We claim that for an appropriate choice of the parameters
$K,L$ the collection of variables $\{Z_k\}$ dominates a product 
measure 
\begin{pro}
\label{prop: domination 1}
There exists $K, L, N_0$ and $\ga >0$ such that for $N \geq N_0$
\begin{eqnarray*}
\forall k \leq M, \qquad
\bbQ \left( Z_k = 1 \big| \; 
Z_j = \eta_j, \quad j \leq k-1 \right) \geq \ga \, , 
\end{eqnarray*}
for any collection of variables $\{\eta_j \}_{j \leq M}$ taking values in $\{0,1\}^M$.  
\end{pro}
The proof is postponed to Section \ref{sec: proof}.

		\subsection{Wired boundary conditions}
		\label{subsec: wired}

Following the previous Subsection, we are going to define
another type of random variables which are
related to the wired FK measure.  
The FK counterpart of the Gibbs measure $\mu^h_{\gb, \gL_N}$
with boundary magnetic field $h>0$ is denoted by 
$\FKm{s,{\rm w}}{\gb,\gL_N}$ and is defined as the wired FK measure
in $\Ew{\gL_N}$ for which a bond $(i,j)$  in $\Ew{\gL_N} \setminus \Ef{\gL_N}$ 
has intensity $s_{(i,j)} = 1 - \exp(-2h J(i-j))$ instead  of $p_{(i,j)}$. 
The intensities of the bonds in $\Ef{\gL_N}$ remain as defined 
in Subsection \ref{subsec: FK}.

\medskip

Using the notation of Definition \ref{def: Z},
we introduce  new random variables indexed by the 
set $\Xi_{N,L}$.
\begin{defi}
\label{def: tilde Z}
For any $x$ in $\Xi_{N,L}$, we declare that $\widehat Z_x =1$ if 
there exists at least one open bond in $\Ew{\gL_N} \setminus \Ef{\gL_N}$ 
joining $T_L(x)$ to $\gL_N$. Otherwise we set $\widehat Z_x =0$.

\medskip

Let $\widehat \bbQ$ be the image measure on $\{0,1\}^{\Xi_{N,L}}$ of 
$\FKm{s,{\rm w}}{\gb,\gL_N}$ by the application $\go \to \{ \widehat Z_x (\go)\}$.
\end{defi}
As in the previous Subsection, the random variables 
$\{\widehat Z_k = \widehat Z(x_k)\}_{k \leq M}$ are ordered wrt the 
lexicographic order in $\Xi_{N,L}$.\\

To any bond configuration $\go$ in $\Ew{\gL_N} \setminus \Ef{\gL_N}$,
one associates two types of boundary conditions: 
$\pi^{\widehat Z(\go)}$ which is defined as in \eqref{eq: pi}  and 
\begin{equation}
\forall b \not \in \Ef{\gL_N}, \qquad
\pi^\go_b = 
\begin{cases}
\go_b, \quad \text{if} \ b \in \Ew{\gL_N} \setminus \Ef{\gL_N} \, ,\\
1, \quad \text{otherwise} \, .
\end{cases}
\end{equation}
Thus the following FKG domination holds $\pi^{\widehat Z(\go)}
\succ  \pi^\go$ and conditionally to the bond configuration outside 
$\Ef{\gL_N}$
\begin{eqnarray*}
\Psi(\widehat Z(\go)) \geq 
\FKm{\pi^\go}{\gb,\gL_N} \big( 0 \lra \partial \gL_N \big) \, ,
\end{eqnarray*}
where $\Psi$ was introduced in \eqref{eq: FK Z}.
This leads to 
\begin{eqnarray}
\label{eq: tilde > omega}
\widehat \bbQ  \big( \Psi(\widehat Z)  \big) \geq 
\FKm{s,{\rm w}}{\gb,\gL_N} \big( 0 \lra \partial \gL_N \big) \, .
\end{eqnarray}

\medskip

Finally, we check that uniformly in $N$ the variables 
$\{\widehat Z_k\}$ satisfy
\begin{pro}
\label{prop: domination 2}
For any collection of variables $\{\eta_j \}_{j \leq M}$ taking values
in $\{0,1\}^M$
\begin{eqnarray*}
\forall k \leq M, \qquad
\widehat \bbQ   \left( \widehat Z_k = 1 \big| \; 
\widehat Z_j = \eta_j ,  \quad j \leq k-1 \right) \leq R L^{d-1} s_h \, , 
\end{eqnarray*}
where $s_h = \max s_{(i,j)}$ and $R$ is the interaction range.
\end{pro}

\begin{proof}
For a given $k \leq M$, the variable $\widehat Z_k$ is an increasing
function supported only by the set of bonds joining $T_L(x_k)$ to $\gL_N$ 
which we denote by $\cT_k$.
From FKG inequality, we have
\begin{eqnarray*}
\widehat \bbQ  \left( \widehat Z_k = 1 \big| \; 
\widehat Z_j = \eta_j \quad j \leq k-1 \right) &\leq& 
\FKm{s,{\rm w}}{\gb,\cT_k} \left( \widehat Z_k (\go) = 1 \right)\\ 
&\leq& 
\FKm{s,{\rm w}}{\gb,\cT_k} \left( \exists \ \text{an open bond in $\cT_k$}
\right) \, . 
\end{eqnarray*}
After conditioning, the $R L^{d-1}$ bonds in $\cT_k$ are independent and
open with intensity at most $s_h$. Thus the Proposition follows. 
\end{proof}

		\subsection{The coupling measure}
		\label{subsec: coupling}

We are going to define a joint measure $\bbP$ for the variables 
$\{ Z_k, \widehat Z_k \}_{k \leq M}$. 
The coupling will be such that 
\begin{eqnarray}
\label{eq: condition 1}
\bbP \ a.s. \quad
\{ Z_k \} \succ \{\widehat Z_k \}, \quad \text{i.e.} \qquad
\bbP \left( \big\{ Z_k \geq \widehat Z_k, \quad \forall k \leq M 
\big\} \right) = 1 \, , 
\end{eqnarray}
and the marginals coincide with $\bbQ$ and $\widehat \bbQ$, i.e. 
for any function $\phi$ in $\{0,1\}^{\Xi_{N,L}}$
\begin{eqnarray}
\label{eq: condition 2}
\bbP \big( \phi(Z) \big) = \bbQ \big( \phi(Z) \big)
\qquad \text{and} \qquad
\bbP \big( \phi(\widehat Z) \big) = 
\widehat \bbQ \big( \phi(\widehat Z) \big) \, .
\end{eqnarray}

\medskip

\begin{pro}
\label{prop: coupling}
There exists $K, L$ and $h>0$ such that for any $N$ large enough, one can
find a coupling $\bbP$ satisfying the conditions \eqref{eq: condition 1} 
and \eqref{eq: condition 2}.
\end{pro}

\begin{proof}

The existence of the coupling is standard and follows from Propositions
\ref{prop: domination 1} and \ref{prop: domination 2}.
First choose $K, L$ large enough such that Proposition \ref{prop: domination 1}
holds and then fix $h$ such that $\ga > R L^{d-1}  s_h$.
The coupling $\bbP$ is defined recursively. 
Suppose that the first $k \leq M-1$ variables
$\cZ_k = \{Z_i\}_{i \leq k}, \widehat \cZ_k = \{\widehat Z_i\}_{i \leq k}$
are fixed such that
$$ \forall i\leq k, \qquad  Z_i \geq \widehat Z_i \, . $$
We define
\begin{eqnarray*}
\left\lbrace
	\begin{array}{l}
\bbP \big( Z_{k+1} = 1, \widehat Z_{k+1} = 0  \; \big| \cZ_k, \widehat \cZ_k
\big) = \bbQ \left( Z_{k+1} = 1 \; \big| \cZ_k \right)
- \widehat \bbQ \big( \widehat Z_{k+1} = 1  \; \big| \widehat \cZ_k \big)
\, , \\
\bbP \big( Z_{k+1} = 1, \widehat Z_{k+1} = 1  \; \big| \cZ_k, \widehat \cZ_k
\big)= \widehat \bbQ \big(\widehat Z_{k+1} = 1  \; \big| \widehat \cZ_k
\big) \, , \\
\bbP \big(  Z_{k+1} = 0, \widehat Z_{k+1} = 0  \; \big| \cZ_k, \widehat \cZ_k
\big)  =   \bbQ \big( Z_{k+1} = 0  \; \big| \widehat \cZ_k \big)
\, .
	\end{array}
\right.
\end{eqnarray*}

Thanks to Propositions \ref{prop: domination 1} and \ref{prop: domination 2}
the measure is well defined and one can check that the conditions 
\eqref{eq: condition 1} and \eqref{eq: condition 2} are fulfilled.
\end{proof}

		\subsection{Conclusion}
		\label{subsec: proof main thm}

For $\gb<\gb_c$ Theorem \ref{thm: main} holds (see Theorem 5.3 (a) in \cite{G2}), 
thus we focus on the case $\gb>\gb_c$.
As the wired FK measure dominates the free FK measure in the FKG sense,
it is enough to prove
\begin{eqnarray}
\label{eq: inegalite}
\FKm{{\rm f}}{\gb} \big( \{ 0 \lra \infty \} \big)
\geq \FKm{{\rm w}}{\gb} \big( \{ 0 \lra \infty \} \big) \, .
\end{eqnarray}

Let us first fix $K, L,h$ such that Proposition \ref{prop: coupling}
holds.
From \eqref{eq: omega > Z} and \eqref{eq: condition 2} 
\begin{eqnarray*}
\FKm{{\rm f}}{\gb} \big( \{ 0 \lra \infty \} \big) \geq
\bbQ \left(\Psi \big( Z \big) \right)
= \bbP \left(\Psi \big( Z \big) \right) \, .
\end{eqnarray*}
As $\Psi$ is an increasing function, we get from \eqref{eq: condition 1}
\begin{eqnarray*}
\bbP \left(\Psi \big( Z \big) \right) \geq 
\bbP \left(\Psi \big( \widehat Z \big) \right) \, .
\end{eqnarray*}
Finally from \eqref{eq: condition 2} and \eqref{eq: tilde >
omega} we conclude that
\begin{eqnarray*}
\bbP \left(\Psi \big( \widehat Z \big) \right) =
\widehat \bbQ \big( \Psi(\widehat Z)  \big) \geq 
\FKm{s,{\rm w}}{\gb,\gL_N} \big( 0 \lra \partial \gL_N \big) \, .
\end{eqnarray*}
Thus  the previous inequalities imply that for any $N$ large
enough
\begin{eqnarray*}
\FKm{{\rm f}}{\gb} \big( \{ 0 \lra \infty \} \big)
\geq \FKm{s,{\rm w}}{\gb,\gL_N} \big( 0 \lra \partial \gL_N \big)
= \mu^{h}_{\gb,\gL_N} (\gs_0)  \, ,
\end{eqnarray*}
where $\mu^{h}_{\gb,\gL_N}$ denotes the Gibbs measure with boundary
magnetic field $h = - \frac{1}{2} \log(1-s)$.
It was proven by Lebowitz \cite{Leb2} 
and Messager, Miracle Sole, Pfister \cite{MMP} that for any $h>0$
\begin{eqnarray*}
\lim_{N \to \infty} \mu^{h}_{\gb,\gL_N} (\gs_0) = \mu^+_\gb (\gs_0) \, .
\end{eqnarray*}
Therefore the correspondence between the Ising model and the FK
representation \eqref{eq: averaged magnetization}
completes the derivation of inequality  \eqref{eq: inegalite}.

%Let $\cC$ be the event that $i$ and $j$ are connected by an open 
%path which does not use the bond $(i,j)$. We write
%\begin{eqnarray*}
%\FKm{{\rm w}}{\gb} \big( \go_{(i,j)} \big) =
%\FKm{{\rm w}}{\gb} \big( \go_{(i,j)} 1_\cC \big)
%+ \FKm{{\rm w}}{\gb} \big( \go_{(i,j)} (1- 1_\cC) \big) \, .
%\end{eqnarray*} 
%Conditionally to $\cC$ or its complement, the distribution 
%of $\go_{(i,j)}$ is determined thus
%\begin{eqnarray*}
%\FKm{{\rm w}}{\gb} \big( \go_{(i,j)} \big) =
%p \FKm{{\rm w}}{\gb} \big( 1_\cC \big)
%+ \frac{p}{\gb + (1 -p) q}
%\left( 1 - \FKm{{\rm w}}{\gb} \big( 1_\cC) \big) \right) \, .
%\end{eqnarray*} 
%On the other hand, the spin correlation can also be expressed
%in terms of the expectation of $\cC$.
%\begin{eqnarray*}
%\mu^+_\gb (\gs_i \, \gs_j)
%= \FKm{{\rm w}}{\gb} (\cC) +
%\FKm{{\rm w}}{\gb} \big( \go_{(i,j)}(1- 1_\cC) \big) 
%=
%\FKm{{\rm w}}{\gb} ( \cC ) + 
%\frac{p}{\gb + (1 -p) q}
%\left( 1 - \FKm{{\rm w}}{\gb} \big( 1_\cC) \big) \right) \, .
%\end{eqnarray*} 

%Since $\gb \to \mu^+_\gb (\gs_i \, \gs_j)$ is continuous, the
%previous equalities imply that $p \to \FKm{{\rm w}}{\gb} \big( \go_{(i,j)}
%\big)$ is also continuous.

	\section{Proof of Proposition \ref{prop: domination 1}}
	\label{sec: proof}

For any $k$, we write $Z_k = Z_{x_k} = X_k Y_k$, where the
random variables $X_k$ and $Y_k$ are defined as follows
\begin{itemize}
\item $X_k =1$ if and only if the conditions (1) and (2) of Definition 
\ref{def: Z} are both satisfied. Otherwise $X_k =0$. 
\item  $Y_k =1$ if and only if the condition (3) of Definition 
\ref{def: Z} is satisfied. Otherwise $Y_k =0$. 
\end{itemize}

For any collection of variables $\{\eta_j \}_{j \leq M}$ taking values
in $\{0,1\}^M$, we set
$$
\cC = \left\{
Z_j = \eta_j, \quad j \leq k-1\right\} \, . 
$$
We are going to prove that for $K,L$ large enough there exists $c_1,c_2 \in [0,1[$ 
(depending on $K,L$) such that
\begin{eqnarray}
\label{eq: X}
&&\bbQ \left( X_k = 0 \big| \; \cC \right)
\leq c_1 \, ,
\\
\label{eq: Y}
&& \bbQ \left( X_k = 1, Y_k =0 \big| \; \cC \right)
\leq c_2  \bbQ \left( X_k = 1 \big| \; \cC \right) \, .
\end{eqnarray}

\medskip

Proposition \ref{prop: domination 1} is a direct consequence of the
previous inequalities. 
First we write
\begin{eqnarray*}
\bbQ \left( Z_k = 0 \big| \; \cC \right)
= \bbQ \left( X_k = 0 \big| \; \cC \right) +
\bbQ \left( X_k = 1, Y_k =0 \big| \; \cC \right)\, .
\end{eqnarray*}
Using \eqref{eq: Y} and \eqref{eq: X}
\begin{eqnarray*}
\bbQ \left( Z_k = 0 \big| \; \cC \right)
\leq  1- (1- c_2) \bbQ \left( X_k = 1 \big| \; \cC \right)
\leq 1- (1-c_2) (1-c_1) \, .
\end{eqnarray*}
Thus  for $K,L$ large enough there is $\ga >0$ such that 
\begin{eqnarray*}
\bbQ \left( Z_k = 1 \big| \; \cC \right)
\geq \ga \, .
\end{eqnarray*}
\qed

\vskip.5cm

\noindent
{\it Proof of \eqref{eq: X}.}

The counterpart for $x_k$ of the site $y$ in Definition \ref{def: Z}
is denoted by $y_k$.
The event  $X_k = 1$ requires first of all that 
\begin{itemize}
\item All the edges in $\bbE \setminus \Ef{\gL_N}$ intersecting $T_L(x_k)$ are open.
\item The $3K/4$ edges $\big\{ \big(x_k + i \vec{n}, x_k + (i+1) \vec{n}
\big) \big\}_{0 \leq i \leq 3K/4}$ are open, where $\vec{n}$ denotes the outward normal to $\gL_{N+1}$ at $x_k$.
\end{itemize}
Let $\cA$ be the intersection of both events.
The support of $\cA$ is disjoint from the support of $\cC$, so that $\cA$ can be satisfied
with a positive probability depending on $K$ and $L$ but not on $\cC$ or $N$.

It remains to check that conditionally to $\cA \cap \cC$, the block 
$\bbB_K(y)$ is good with a positive 
probability depending on $K$. We stress the fact that this statement is not a
direct consequence of \eqref{eq: Peierls}  because $\cA$ cannot be expressed in 
terms of the coarse grained variables. Nevertheless $\cA$ is increasing, thus
one can use similar arguments as in Theorem 3.1 of \cite{pisztora} to conclude
that the estimate \eqref{eq: Peierls} remains valid despite the
conditioning by $\cA$.

Combining the previous statements, we deduce that  \eqref{eq: X} holds with a 
constant $c_1 < 1$.

\vskip.5cm

\noindent
{\it Proof of \eqref{eq: Y}.}

\begin{figure}[h]
\begin{center}
\leavevmode
\epsfysize = 5 cm
\psfrag{gL}[r]{$\gL_{N + 3K/2}^c$}
\psfrag{L}[T]{$L$}
\psfrag{y1}[Tl]{$\bbB_K(y_1)$}
\psfrag{y2}[Tr]{$\bbB_K(y_2)$}
\psfrag{yk}[T]{$\bbB_K(y_k)$}
\psfrag{C}[Tr]{$\frC$}
\psfrag{P}[B]{$\partial \gga$}
\psfrag{ga}[T]{$\gga$}
\epsfbox{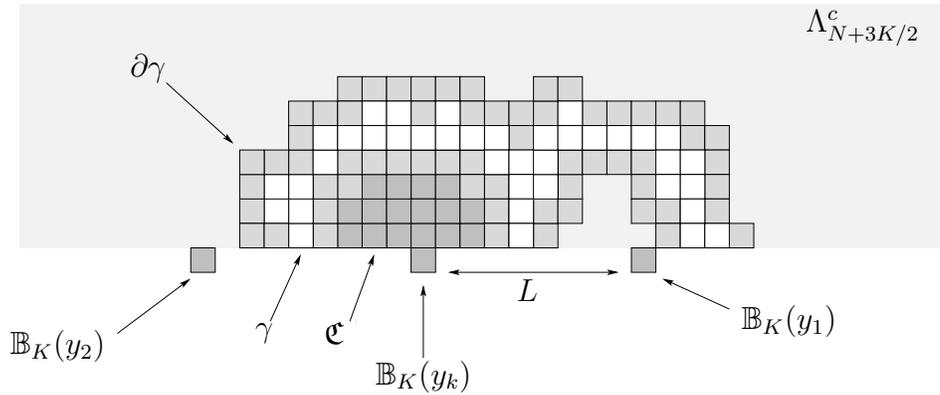}
\end{center}
\caption{The support of the contour $\gG$ is $\gga \cup \partial \gga$
and is included in $\gL_{N + 3K/2}^c$ (the light gray region). 
The blocks $\bbB_K(y_k)$ and $\bbB_K(y_1)$ are disconnected from infinity by $\gG$.
The event $Y_2= 1$ associated to the block $\bbB_K(y_2)$ 
is not determined by $\gG$.}
\label{fig: contour}
\end{figure}

Let $y_k$ be the counterpart of the site $y$ in Definition \ref{def: Z}.
If $Y_k =0$ then there exists $\gG$ a contour of bad blocks in
$\gL_{N+ 3K/2}^c$ disconnecting $y_k$ from infinity 
(see (3) of Definition \ref{def: Z}). 
More precisely, we define the contour $\gG$ as follows.
Let $\frC$ be the maximal connected component of  
good blocks in $\gL_{N+ 3K/2}^c$ connected to $\bbB_K(y_k)$.
If $Y_k =0$, $\frC$ is finite and $\gga$ is defined as the support of the maximal 
$\star$-connected component of bad blocks in $\gL_{N+ 3K/2}^c$
which intersects the boundary of $\frC$ or simply the block connected to
$\bbB_K(y_k)$ if $\frC$ is empty.
By construction the boundary of $\gga$, denoted by $\partial \gga$, contains only 
good blocks. The contour $\gG$ is defined as the intersection of the
events $\gG_0$ and $\gG_1$, where the configurations in  $\gG_0$
contain only bad blocks in $\gga$ and those in  $\gG_1$ contain only
good blocks in $\partial \gga$ (see figure \ref{fig: contour}).

We write
\begin{eqnarray}
\label{eq: Y 1} 
\bbQ \left( \{X_k = 1 \} \cap \{Y_k =0\} \cap \cC \right)
\leq \sum_{\gG}
\FKm{{\rm f}}{\gb} \left(\gG \cap \{ X_k = 1 \} \cap \cC \right) \, ,
\end{eqnarray}
where the sum is over the contours in $\gL_{N+3 K/2}^c$ surrounding $y_k$.

For a given $\gG$, we are going to prove
\begin{eqnarray}
\label{eq: Y 2} 
\FKm{{\rm f}}{\gb} \left(\gG \cap \{X_k = 1 \} \cap \cC \right) 
\leq \exp \left( - \frac{C}{2} | \gG | \right)
\FKm{{\rm f}}{\gb} \left( \{X_k = 1\} \cap \cC \right) \, ,
\end{eqnarray}
where $C = C(K,\gb)$ was introduced in \eqref{eq: Peierls}  and $|\gG|$ stands
for the number of blocks in $\gga$.
For $K$ large enough, the constant $C$ can be chosen arbitrarily large
so that the combinatorial factor arising by summing over the contours $\gG$
in \eqref{eq: Y 1} remains under control.
This implies that there exists $c_2 \in ]0,1[$ such that
\begin{eqnarray*}
\bbQ \left( \{X_k = 1 \} \cap \{Y_k =0\} \cap \cC \right)
\leq c_2 \; \FKm{{\rm f}}{\gb} \left( \{ X_k = 1 \} \cap \cC
\right) 
\, .
\end{eqnarray*}
Thus the inequality  \eqref{eq: Y} follows.

\medskip

In order to prove \eqref{eq: Y 2}, we specify the set $\cC$ and
for notational simplicity assume that it is of the form
$\cC =\cC_0 \cap \cC_1$ with 
$$
\cC_0 = \left\{Z_j = 0, \quad j \leq k_0 \right\},
\qquad
\cC_1 = \left\{Z_j = 1, \quad k_0+1 \leq j \leq k-1 \right\}\, . $$
The difficulty  to derive \eqref{eq: Y 2} is that $\gG$ may
contribute to the event $\cC_0$
so that a Peierls argument cannot be applied directly. 
For this reason we decompose $\cC_0$ into $2^{k_0}$ disjoint sets for
which the state of the first $k_0$ variables is prescribed such that
either $\{X_j = 1, Y_j=0\}$ or $\{X_j = 0\}$.
Once again for simplicity we will only consider the subset
$\cD = \cD_0 \cap \cD_1$ of $\cC_0$ such that 
$$
\cD_0 = \left\{X_j = 1, Y_j=0, \quad j \leq k_1 \right\},
\qquad
\cD_1 = \left\{X_j = 0, \quad k_1 +1 \leq j \leq k_0 \right\}\, . 
$$
The derivation of \eqref{eq: Y 2} boils down to prove the estimate
below
\begin{eqnarray}
\label{eq: Y 3} 
\FKm{{\rm f}}{\gb} \left( 
 \gG \cap \{X_k = 1\} \cap \cD \cap \cC_1 \right) \leq 
\exp \left( - \frac{C}{2} | \gG | \right)
\FKm{{\rm f}}{\gb} \left( \{X_k = 1\} \cap \cD \cap  \cC_1 \right) 
\, .
\end{eqnarray}
Finally, we suppose that $\cD_0$ is such that the first $k_2$ sites
$\{ y_j \}_{j \leq k_2}$ are disconnected from infinity by $\gG$
and the others $k_1-k_2$ are not surrounded by $\gG$ (see figure 
\ref{fig: contour}).
Notice that erasing the contour $\gG$ may affect the state of the first
$k_2$ sites, but not of the other $k_1 - k_2$.
By construction, if $\cE = \left\{X_j = 1, Y_j=0, \quad k_2 +1 \leq j 
\leq k_1 \right\}$, then 
\begin{eqnarray*}
\FKm{{\rm f}}{\gb} \left(\gG \cap \{X_k = 1\} \cap \cD \cap \cC_1 \right) 
=
\FKm{{\rm f}}{\gb} \left( \gG \cap \{X_k = 1\} \cap
\left\{X_j = 1, \ \  j \leq k_2 \right\}
\cap \cE \cap \cD_1 \cap  \cC_1 \right) 
\, .
\end{eqnarray*}
Conditionally to $\gG_1$,
all the events in the RHS are independent of $\gG_0$ so that by
conditioning wrt the configurations in $\partial \gga$, one can apply
the Peierls bound \eqref{eq: Peierls} 
\begin{eqnarray*}
&& \FKm{{\rm f}}{\gb} \left(\gG \cap \{X_k = 1\} \cap \cD \cap \cC_1 \right)\\
&& \qquad  \leq 
\exp  \left( - C | \gG | \right) \;
\FKm{{\rm f}}{\gb} \left(\gG_1 \cap  \{X_k = 1\} \cap
\left\{X_j = 1, \quad j \leq k_2 \right\}
\cap \cE \cap \cD_1 \cap  \cC_1 \right) 
\, .
\end{eqnarray*}
By modifying the bonds around each block $\bbB_K(y_j)$
one can recreate the events $\{Y_j =0\}_{j \leq k_2}$ and thus
$\cD$. 
First of all notice that $\gG_1$ screens the blocks  $\bbB_K(y_j)$
from the other events in the RHS. 
Thus one can turn the blocks in $\gL^c_{N + 3K/2}$ connected to 
each site $\{y_j\}_{j \leq k_2}$ into bad blocks without 
affecting the event below 
$$ \{X_k = 1\} \cap
\left\{X_j = 1, \quad j \leq k_2 \right\}
\cap \cE \cap \cD_1 \cap  \cC_1.$$
%Notice that this might not preserve the event $\gG_1$.
For each block, this has a cost $\ga_K$ depending only on K (and $\gb$)
\begin{eqnarray*}
\FKm{{\rm f}}{\gb} \left(\gG \cap \{X_k = 1\} \cap \cD \cap \cC_1 \right)
\leq 
\exp \left( - C | \gG | \right) \; \big( \ga_K \big)^{k_2} \;
\FKm{{\rm f}}{\gb} \left(  \{X_k = 1\} \cap \cD \cap  \cC_1 \right) 
\, .
\end{eqnarray*}
By construction, the distance between each site $\{y_j\}_{j \leq k_2}$
is at least $L = \ell K$. The contour $\gG$ surrounds  $k_2$ sites
in $\Xi_{N,L}$ so that $|\gG|$ must be larger than $\ell k_2$
(see figure \ref{fig: contour}).  
Therefore for $\ell$ large enough, the Peierls bound compensates the 
cost  $\big( \ga_K \big)^{k_2}$
\begin{eqnarray*}
\FKm{{\rm f}}{\gb} \left(\gG \cap \{X_k = 1\} \cap \cD \cap \cC_1 \right)
\leq 
\exp \left( - \frac{C}{2} | \gG | \right) \;
\FKm{{\rm f}}{\gb} \left( \{X_k = 1\} \cap \cD \cap  \cC_1 \right) 
\, .
\end{eqnarray*}
This completes  \eqref{eq: Y 3}.
Similar results would be valid for any decomposition of the set $\cC$.
In particular $\cC_0$ can be represented as the disjoint union of the type
$\cC_0 =  \bigvee_{\cD_0,\cD_1} \cD_0 \cap \cD_1$, thus summing over 
the sets $\cD$, we derive  \eqref{eq: Y 2}.

\end{document}